\documentclass[12pt]{amsart}
\usepackage{latexsym,enumerate}
\usepackage{amssymb, xcolor,mathrsfs}
\usepackage{amsmath,amsthm,amsfonts,amssymb,latexsym, mathabx}

\newtheorem{theorem}{Theorem}
\newtheorem{corollary}{Corollary}
\newtheorem{lemma}{Lemma}

\newtheorem{conjecture}{Conjecture}

\newtheorem*{mainthm}{Main Theorem}

\theoremstyle{definition}

\newtheorem*{xrem}{Remark}

\newcommand{\Z}{\mathbb{Z}}
\newcommand{\Q}{\mathbb{Q}}
\newcommand{\R}{\mathbb{R}}

\newcommand{\la}{\langle}
\newcommand{\ra}{\rangle}

\newcommand{\m}{\;\mbox{mod}\;}

\DeclareMathOperator{\lcm}{lcm}

\DeclareMathOperator{\Gal}{Gal}

\begin{document}
	
	\baselineskip=17pt
	
	%%%%%%%%%%%%%%%%

	\title[Answer to a Question of Hung And Tiep]{Answer to a Question of Hung and Tiep on Conductors of Cyclotomic Integers}
	
	\author[Christopher Herbig]{Christopher Herbig}
	\address{%
		Northern Illinois University\\
		Department of Mathematical Sciences\\
		Dekalb, IL 60115\\
		USA}
	\email{cherbig@niu.edu}
	\date{}
	
	\maketitle
	
	%% Classification and key words; note that the 2010 classification is used:
	
	\renewcommand{\thefootnote}{}
	
	\footnote{2020 \emph{Mathematics Subject Classification}: Primary 11R18; Secondary 20C15.}
	
	\footnote{\emph{Key words and phrases}: cyclotomic extensions, roots of unity, character values.}
	
	\renewcommand{\thefootnote}{\arabic{footnote}}
	\setcounter{footnote}{0}
	
	%%%%%%%%
	
	\begin{abstract}
		In Question 5.2 of \cite{HungTiep}, Hung and Tiep asked the following: If $\alpha$ is a sum of $k$ complex roots of unity and $\Q_{c(\alpha)}$ is the smallest cyclotomic field containing $\alpha$, is it true that $|\Q_{c(\alpha)}:\Q(\alpha)| \leq k$? We answer this question in the negative. Using known results on minimal vanishing sums, we also characterize all cyclotomic integers with $k \leq 4$ for which the inequality fails. In \S 6, we bound the growth of $|\Q_{c(\alpha)}:\Q(\alpha)|$ as a function of $k$.
	\end{abstract}
	
	\section{Introduction}
	
	We discuss a question of Hung and Tiep regarding the relationship between weights and conductors of cyclotomic integers posed in Question 5.2 of \cite{HungTiep}. Although the question was asked in the context of group characters and their fields of values, this question may be reformulated into a conjecture which is purely number-theoretic.
	
	Before proceeding, we shall introduce some necessary notation. We let $\Q_n$ denote the extension of $\Q$ generated by a primitive $n$th root of unity. For a complex character $\chi$ of a finite group, we define $\Q(\chi)$ to be the field generated by all the values that $\chi$ takes on the group. By a \emph{cyclotomic integer}, we mean any algebraic integer which may be expressed as a sum of roots of unity. For a particular sum of roots of unity, say $\alpha = \sum_{i=1}^k\epsilon_i$, we define the \emph{weight} of $\alpha$, denoted $w(\alpha)$, to be the integer $k$. We define the \emph{length} of a cyclotomic integer $\alpha$, denoted $l(\alpha)$, to be the smallest possible weight among all representations of $\alpha$ as a sum of roots of unity. By the \emph{conductor} of a cyclotomic integer $\alpha$, we mean the smallest possible integer $c(\alpha)$ such that $\Q(\alpha) \subseteq \Q_{c(\alpha)}$. We similarly define $c(\chi)$ for a complex character $\chi$ to mean the smallest possible integer such that $\Q(\chi) \subseteq \Q_{c(\chi)}$. Lastly, we will use $\zeta_m$ to denote the complex root of unity $e^{\frac{2\pi i}{m}}$.
	
	Hung and Tiep's Question 5.2 may be stated as follows:
	
	\begin{conjecture}\label{C:1}
		\textup{\cite[Question 5.2]{HungTiep}} Let $G$ be a finite group, $\chi$ an irreducible complex character of $G$, and $g \in G$. Then, we have
		
		$$
		|\Q_{c(\chi(g))}: \Q(\chi(g))| \leq \chi(1).
		$$
	\end{conjecture}

	Hung and Tiep introduced the above conjecture as the element-level version of the below conjecture, which deals with conductors of the field generated by all values of a group character:
	
	\begin{conjecture}\label{C:2}
		\textup{\cite[Conjecture 1.2]{HungTiep}} Let $G$ be a finite group and $\chi$ an irreducible complex character of $G$. Then, 
		
		$$
		|\Q_{c(\chi)} : \Q(\chi)| \leq \chi(1).
		$$
	\end{conjecture}

	In \cite{HungTiep}, Hung and Tiep have verified Conjecture 2 for the cases where $G$ was alternating, unitary, or a general linear group. They also proved that Conjecture \ref{C:2} holds whenever $\chi(1) \leq 3$. More recently, Hung, Tiep and Zalesski have verified the conjecture in \cite{HungTiepZal} when $\chi(1) = p$ for some prime $p$.
	
	Conjecture 1 may be equivalently reformulated into a purely number-theoretic conjecture.
	
	\begin{conjecture}\label{C:3}
		If $\alpha$ is a cyclotomic integer, then
		
		$$[\Q_{c(\alpha)}:\Q(\alpha)] \leq l(\alpha).$$
	\end{conjecture} 
	
	The equivalence of Conjectures \ref{C:1} and \ref{C:3} is seen by showing that for every sum of roots of unity, there exists an irreducible character of some group taking on such a value. This is done by considering the irreducible characters of the imprimitive regular wreath product $G:= C_m \wr C_n$. We choose a generator $\lambda$ of the character group of $C_m$ and set $\psi = \lambda \times \prod_{i=2}^n 1_{C_m}$, which we view as an irreducible character of the base group $\prod_{i=1}^n C_m$. One sees that the inertia group $I_G(\psi)$ is merely the base group, and so, the induced character $\psi^G$ is an irreducible character of the group $G$. One then verifies the values of $\psi^G$ range over all sums of $m$th roots of unity having weight $n$.
	
	Since our proofs require no character theory, we will exclusively treat Conjecture \ref{C:3} for the remainder of this paper. We will produce counterexamples for all lengths at least 3 and show the conjecture holds when the length is 2. We will then classify the counterexamples having length 3 or 4. This will be accomplished by using the classification of minimal vanishing sums of roots of unity having weight at most 8. This may be found in \cite{ChrisDyk}, \cite{ConJon}, or \cite{PoonRub}. Although Conjecture \ref{C:3} turns out to be false, it is natural to wonder about the asymptotic behavior of $d(k) := \max\{ |\Q_{c(\alpha)}: \Q(\alpha)| \; \big| \; l(\alpha) = k\}$. It is not immediately clear that $d(k)$ is finite. However, using our methods, we are able to bound $d(k)$ by a function of $k$.
	
	\begin{mainthm}
		For $k \in \Z^+$, we have	
		$$d(k) < k!\bigg( \dfrac{2.376ke^{4ck}}{\log(k+1)}\bigg)^k,$$
		where $c = 1.000028$. In particular, $d(k)$ is finite for all $k$.
	\end{mainthm}
	
	The proof will be given in \S 6, and we comment on the bound in the conclusion.
	
	\section{Cyclotomic Integers Of Weight At Least 5}
	
	A straightforward counterexample to Conjecture \ref{C:3} is found by taking $$\alpha = \zeta_8 + \zeta_8^7 + \zeta_7 + \zeta_7^2 + \zeta_7^4.$$ One establishes $c(\alpha) = 56$ by using the fact that $\Q_n \cap \Q_m = \Q_{\gcd(n,m)}$. Indeed, if $\alpha \in \Q_{28}$, this implies $\zeta_8 + \zeta_8^7 = \sqrt{2}$ would lie in $\Q_8 \cap \Q_{28} = \Q(i)$, a contradiction. Similarly, $\alpha \notin \Q_8$. This suffices to see that $\alpha$ does not lie in any proper cyclotomic subfield of $\Q_{56}$. We now claim that $[\Q_{56}:\Q(\alpha)] \geq 6$. Define $\sigma_1, \sigma_2 \in \Gal(\Q_{56})$ by $\sigma_1(\zeta_{56}) = \zeta_{56}^9$ and $\sigma_2(\zeta_{56}) = \zeta_{56}^{15}$. One sees that $\sigma_1$ and $\sigma_2$ permute the terms of $\alpha$, and in particular, $\sigma_1 , \sigma_2 \in \Gal(\Q_{56}/\Q(\alpha))$. However, $\la \sigma_1 , \sigma_2 \ra$ is cyclic of order 6, despite the fact $w(\alpha) = 5$.
	
	An infinite family of counterexamples can be produced by similar means. In particular, we prove the following theorem:
	
	\begin{theorem}\label{T:1}
		Let $n$ be a composite, square-free integer which is divisible by  neither 2 nor 3. Then there exists a cyclotomic integer $\alpha$ having conductor $n$ for which $|\Gal(\Q_{n}/\Q(\alpha))| > l(\alpha)$.
	\end{theorem}
	
	\noindent\emph{Proof:} We proceed inductively. Let $p$ and $q$ be distinct primes, both at least 5. Choose nontrivial, proper subgroups $H_p$ and $H_q$ of $\Gal(\Q_p)$ and $\Gal(\Q_q)$ respectively. Assume further that either $H_p$ or $H_q$ has order at least 3. Such a choice is possible by our assumptions on $p$ and $q$. Set
	
	$$
	\alpha_1 = \sum_{\sigma \in H_p} \zeta_p^\sigma, \;\;\; \alpha_2 = \sum_{\sigma \in H_q} \zeta_q^\sigma,
	$$
	and
	$$
	\alpha = \alpha_1 + \alpha_2.
	$$
	
	It is straightforward to verify that $c(\alpha) = pq$, namely, since $H_p$ and $H_q$ were chosen to be proper subgroups, neither $\alpha_1$ nor $\alpha_2$ can be rational. By the Chinese Remainder Theorem, it is possible to choose automorphisms $\tau_p, \tau_q \in \Gal(\Q_{pq})$ such that $\tau_p$ restricts to a generator of $H_p$ in $\Gal(\Q_p)$ and the identity in $\Gal(\Q_q)$ and that $\tau_q$ restricts to the identity in $\Gal(\Q_p)$ and a generator of $\Gal(\Q_q)$. Thus, $\la \tau_p , \tau_q \ra \cong H_p \times H_q$, so that
	$$
	|\Gal(\Q_{pq}/\Q(\alpha))| \geq |H_p \times H_q| > |H_p| + |H_q| = w(\alpha) \geq \ell(\alpha),
	$$
	where the strict inequality holds since either $H_p$ or $H_q$ has order at least 3.
	
	Now, assume that there exists a counterexample $\alpha'$ with $c(\alpha') = n$ where $n$ is a square-free integer having $k$ prime factors. Let $r \geq 5$ be a prime such that $\gcd(n,r) = 1$. Set $H_n = \Gal(\Q_n/\Q(\alpha'))$, and let $H_r$ be a nontrivial, proper subgroup of $\Gal(\Q_r)$. Set
	
	$$
	\alpha = \alpha' + \sum_{\sigma \in H_r} \zeta_r^\sigma.
	$$
	
	Similarly to the base case, we obtain $c(\alpha) = nr$, and it is possible to use the Chinese Remainder Theorem to produce a subgroup of $\Gal(\Q_{nr})$ isomorphic to $H_n \times H_r$. Now,
	
	$$
	|\Gal(\Q_{nr}/\Q(\alpha))| \geq |H_n \times H_r| > |H_n| + |H_r | > w(\alpha') + |H_r| \geq w(\alpha) \geq l(\alpha),
	$$
	which proves the theorem. $\square$\\
	
	\begin{xrem}
		This construction does not extend easily if $n$ above is replaced with an integer which is not square-free, although counterexamples do exist such as the first example we gave. The obstruction comes from the fact that it is not easy in general to determine the conductor of the above orbit sums when the square-free assumption is omitted. This construction also cannot be used to construct counterexamples of weight 3 or 4. Our method only produces automorphisms which permute the terms of an $\alpha$, and since the largest possible orders for abelian subgroups of $S_3$ and $S_4$ are 3 and 4 respectively, Galois automorphisms arising from permutations of the roots alone will not suffice to produce counterexamples.
	\end{xrem}  
	
	\begin{corollary}\label{CL:1}
		For each $k \geq 5$, there exists a sum of roots of unity $\alpha$ such that $w(\alpha) = k$ and Conjecture \ref{C:3} is false.
	\end{corollary}
	
	\noindent\emph{Proof:} Write $k = i + j$ where $i$ is at least 2 and $j$ is at least 3. We are able to produce primes $p$ and $q$ such that $i$ and $j$ properly divide $p - 1$ and $q - 1$ respectively by Dirichlet's Theorem on arithmetic progressions. Using the notation established in Theorem \ref{T:1}, we choose $H_p \subset \Gal(\Q_p)$ of order $i$ and $H_q \subset \Gal(\Q_q)$ of order $j$, and we set
	
	$$
	\alpha = \sum_{\sigma \in H_p} \zeta_p^\sigma + \sum_{\sigma \in H_q} \zeta_q^\sigma.
	$$  Similarly to the base case of Theorem \ref{T:1}, this constitutes a counterexample to Conjecture \ref{C:3}. $\square$\\
	
	We will now use Theorem \ref{T:1} to show $d(k)$ grows at least exponentially. 
	
	\begin{corollary}\label{CL:2}
		Let $k \in \Z$ with $k \geq 6$. There exists a cyclotomic integer $\alpha$ such that $w(\alpha) = k$ and
		
		$$
		|\Q_{c(\alpha)} : \Q(\alpha)| \geq 3^{\lfloor k/3\rfloor}
		$$
		
	\end{corollary}
	\noindent\emph{Proof:} First, assume that $3 \,|\, k$. Choose distinct primes $p_1, p_2, \ldots, p_\ell$ such that $p_i \equiv 1 \m 3$ for each $i$. Choose $H_i \subset \Gal(\Q_{p_i})$ such that $|H_i| = 3$ for each $i$. Now, set
	
	$$
	\alpha = \sum_{i = 1}^{k/3} \sum_{\sigma \in H_i} \zeta_{p_i}^\sigma
	$$
	
	Similarly to the proof of Theorem \ref{T:1}, we may choose Galois automorphisms which permute one of the orbit sums while leaving the remaining orbit sums fixed. Thus, $\Gal(\Q_{c(\alpha)}/\Q(\alpha))$ is a direct product of $k/3$ copies of $C_3$, and so, $|\Q_{c(\alpha)} : \Q(\alpha)| \geq 3^{k/3}$. Otherwise, if $k$ is indivisible by 3, then choose either one or two primitive $q$th roots of unity where $q$ is a prime not in $\{p_i \; | \; i = 1,\ldots , k/3\}$. Add the root(s) of unity chosen to the expression given for $\alpha$ above, and the above inequality will remain true. $\square$

	\section{Some Special Cases}
	
	We now turn to discussing special cases for which the conjecture is true. In particular, if $\alpha$ is a sum of $p^a$th roots of unity or if $l(\alpha) \leq 2$, then the conjecture holds. Before proceeding, we state an important theorem of Lam and Leung on minimal vanishing sums of roots of unity. By a \emph{minimal vanishing sum}, we mean a vanishing sum of roots of unity such that no proper subset of the terms sums to zero. 
	
	\begin{lemma}\label{L:1}
		\textup{\cite[Theorems 4.5 \& 6.5]{LamLe}} Let $\alpha$ be a minimal vanishing sum of roots of unity, and let $\pi$ be the set of all primes dividing the order of any one of the terms. Then one of the following situations occurs:
		\begin{enumerate}
			\item[\textup{(i)}] $\alpha$ is of the form
			\begin{equation}\label{E:1}
				\epsilon\sum_{i=0}^{p-1}\zeta_p^i,			
			\end{equation}
			where $\epsilon$ is some root of unity and $p$ is some prime.
			\item[\textup{(ii)}]
			Otherwise, $|\pi| \geq 3$ and $\alpha$ satisfies
			
			\begin{equation}\label{E:2}
				w(\alpha) \geq (p_1 - 1)(p_2 - 1) + p_3 - 1 \geq p_3 + 1 \geq 6,
			\end{equation}		
			where $p_1 < p_2 < p_3$ are the three smallest elements of $\pi$. Moreover, equality holds in \eqref{E:2} iff
			
			\begin{equation}\label{E:3}
				\alpha = \epsilon\bigg[(\sum_{i = 1}^{p_1 - 1}\zeta_{p_1}^i)(\sum_{i = 1}^{p_2 - 1}\zeta_{p_2}^i) + \sum_{i = 1}^{p_3 - 1}\zeta_{p_3}^i\bigg],
			\end{equation}
			where $\epsilon$ is some root of unity.
		\end{enumerate} 
	\end{lemma} 
	
	We will use this lemma to prove another lemma which will allow us to deal with the case when $\alpha$ is a sum of $p^a$th roots of unity for some positive integer $a$.
	
	\begin{lemma}\label{L:2} Let $\alpha = \sum_{i=1}^k \epsilon_i$ where each $\epsilon_i$ is a (not necessarily primitive) $2^ap^b$th root of unity for some prime $p$ and positive integers $a$ and $b$. Then, $\Gal(\Q_{2^ap^b}/\Q(\alpha))$ permutes the terms of $\alpha$.
	\end{lemma}
	
	\noindent\emph{Proof:} It is of no loss to assume that the sum $\sum_{i=1}^k \epsilon_i$ is chosen such that there exists no sum of $2^ap^b$th roots of unity equal to $\alpha$ having strictly smaller weight. Choose $\sigma \in \Gal(\Q_{2^ap^b}/\Gal(\alpha))$. We have that $$\alpha - \sigma(\alpha) = \sum_{i=1}^k \epsilon_i - \sum_{i=1}^k \epsilon_i^\sigma$$ is a vanishing sum of $2k$ $2^{a}p^b$th roots of unity. We can thus partition the set $\{\epsilon_1 , \epsilon_2 , \ldots , \epsilon_k , -\epsilon_1^\sigma , \ldots , -\epsilon_k^\sigma\}$ into sets of minimal vanishing subsums. By our minimality assumption, each set in the partition consists of the same number of terms of $\alpha$ as of $\sigma(\alpha)$. Indeed, if this were not true, there would be at least one set in the partition containing more terms of $\alpha$ than of $\sigma(\alpha)$, and it would then be possible to replace the terms of $\alpha$ from this set with fewer terms from $\sigma(\alpha)$, contradicting the minimality of $\alpha$. In particular, the minimal vanishing subsums of $\alpha - \sigma(\alpha)$ each have even weight.
	
	Since the orders of the terms consist of at most two prime divisors, case (i) of Lemma \ref{L:1} applies, and the minimal vanishing sums will have weight either 2 or $p$. However, each minimal vanishing subsum has even weight and must then be of the form $\epsilon_i - \epsilon_i$, again by Lemma \ref{L:1}. Since one of the terms of a given minimal vanishing sum is a term of $\alpha$ and the other is a term of $\sigma(\alpha)$, it now follows that $\sigma$ permutes the terms. Since $\sigma$ was arbitrary, the lemma follows. $\square$\\
	
	\begin{xrem} This lemma did not require the full force of Lemma \ref{L:1}, but merely an earlier result of de Bruijn in \cite{Bru} which states that minimal vanishing sums of $p^aq^b$th roots of unity always take the form of \eqref{E:1} in Lemma \ref{L:1}. Notice that we cannot replace $2^ap^b$ in the statement of Lemma \ref{L:2} with $p^aq^b$ since the act of subtracting roots of unity having odd order gives rise to roots of unity having twice their original order.
	\end{xrem}
	
	We are now able to prove the following.
	
	\begin{theorem}\label{T:2} Conjecture \ref{C:3} holds whenever $\alpha$ is a sum of $p^a$th roots of unity for some prime $p$ and some positive integer $a$.
	\end{theorem}
	
	\noindent\emph{Proof:} Write $\alpha = \sum_{i=1}^k \epsilon_i$ where each $\epsilon_i$ is a $p^a$th root of unity, and assume that $\alpha$ cannot be represented as a sum of less than $k$ $p^a$th roots of unity. Also, assume that $\Q_{p^a}$ is the smallest cyclotomic field in which $\alpha$ is contained. In particular, at least one $\epsilon_i$ is a primitive $p^a$th root of unity.  Set $G = \Gal(\Q_{p^a}/\Q(\alpha))$. Now, $G$ permutes the terms of $\alpha$ by Lemma \ref{L:2}. We claim that $G$ is cyclic. If this is true, then since $G$ acts regularly on the primitive $p^a$th roots of unity, the orbit containing any primitive $p^a$th root of unity in the sum will have size $|G|$ which then implies $w(\alpha) \geq |G|$ as desired. It is well known that $G$ is cyclic if $p$ is odd. On the other hand, if $p = 2$, then $G$ is identified with a subgroup of $C_2 \times C_{2^{a-2}}$. $G$ cannot contain the subgroup of order 2 in $1 \times C_{2^{a-2}}$ because this would contradict the fact that $c(\alpha) = p^a$. There are only two subgroups of $C_2 \times C_{2^{a-2}}$ not containing this subgroup, and they are both cyclic of order 2. $\square$\\
	
	\begin{theorem}\label{T:3}
		Conjecture \ref{C:3} holds for cyclotomic integers of length 2.
	\end{theorem} 
	
	\noindent\emph{Proof:} Assume $\alpha$ is a cyclotomic integer of length 2 and $c(\alpha) = n$. Let $\sigma \in \Gal(\Q_n/\Q(\alpha))$ so that $\alpha - \sigma(\alpha)$ is a vanishing sum of four roots of unity. We partition the terms into minimal vanishing sums. Lemma \ref{L:1} implies that there are no minimal vanishing sums of weight 4, so the only possible partition occurs when the minimal vanishing sums have weight 2. Thus, $\Gal(\Q_n/\Q(\alpha))$ either fixes or swaps the terms of $\alpha$. This forces the Galois group to have order at most 2, as desired. $\square$
	
	\section{Cyclotomic Integers Of Length 3}
	
	Conjecture \ref{C:3} is generally false for cyclotomic integers of length 3. Set
	$$\alpha = i(\zeta_6 + \zeta_5 + \zeta_5^4).$$
	One sees that $\alpha$ has conductor equal to $\lcm[4,6,5] = 60$. The Galois automorphism $\sigma$ of $\Q_{60}$ given by $\sigma(\zeta_{60}) = \zeta_{60}^{47}$ has order 4 and simultaneously fixes $\alpha$. The fact that $\sigma$ fixes $\alpha$ follows from the fact that
	
	$$
	\alpha - \sigma(\alpha) = \big[i(\zeta_6 + \zeta_5 + \zeta_5^4)\big]  - \big[-i(\zeta_6^5 + \zeta_5^2 + \zeta_5^3)\big]  = 0.
	$$
	
	Since there is an essentially unique minimal vanishing sum of roots of unity having weight 6, it is possible to characterize completely the length 3 counterexamples to Conjecture \ref{C:3}.

	\begin{theorem}\label{T:4}
		Let $\alpha$ be a cyclotomic integer of length 3. Then, $$[\Q_{c(\alpha)}:\Q(\alpha)] \leq 4$$ with equality holding iff $\alpha$ is of the form
		
		$$
		\epsilon (\zeta_6^c + \zeta_5^a + \zeta_5^{4a}),
		$$
		where $c \in \Z_6^\times$, $a \in \Z_5^\times$, and $\epsilon$ is some root of unity whose order is a multiple of 4 not divisible by either 3 or 5.
		
	\end{theorem}
	
	\noindent\emph{Proof:} Write $\alpha = \epsilon_1 + \epsilon_2 + \epsilon_3$ for roots of unity $\epsilon_i$, and set
	$$G = \Gal(\Q(\epsilon_1, \epsilon_2, \epsilon_3)/\Q(\alpha)).$$ 
	For $\sigma \in G$, $\alpha - \sigma(\alpha)$ is a vanishing sum of six roots of unity. The terms of this sum may be partitioned into minimal vanishing sums. If the sum is partitioned into three minimal vanishing sums of size 2 for each $\sigma \in G$, then, as in Theorem \ref{T:3}, $G$ permutes the terms of $\alpha$. It follows that $G$ is an abelian subgroup of $S_3$ and therefore has order at most 3. We may then assume that there exists some $\sigma \in G$ that does not permute the terms of $\alpha$. 
	
	The vanishing sum may split into two minimal vanishing sums of size 3. In this case, one of the two minimal vanishing sums must share two terms with $\alpha$. However, we may replace these two terms of $\alpha$ with a single term from $\sigma(\alpha)$, contradicting our assumption that $l(\alpha) = 3$. Also, by Lemma \ref{L:1}, there are no minimal vanishing sums of weight 4, so a partition of $\alpha - \sigma(\alpha)$ into sets of size 4 and 2 is impossible.
	
	The last possible case is where there exists a $\sigma \in G$ such that $\alpha - \sigma(\alpha)$ is itself a minimal vanishing sum. Since this sum has weight 6, Lemma \ref{L:1} implies that the sum will take the form of \eqref{E:3} with $p_1$, $p_2$, and $p_3$ being 2, 3, and 5 respectively. So, we obtain
	
	\begin{equation}\label{E:4}
		\alpha - \sigma(\alpha) = \epsilon(\zeta_6 + \zeta_6^5 + \zeta_5 + \zeta_5^2 + \zeta_5^3 + \zeta_5^4)
	\end{equation}
	where $\epsilon$ is some root of unity. This minimal vanishing sum is unique up to the choice of $\epsilon$, so it follows that the terms of $\alpha$ consist of 3 of the terms of \eqref{E:4}. Using the fact that $l(\alpha) = 3$, one verifies that the terms of $\alpha/\epsilon$ consist of exactly one primitive sixth root of unity and two primitive fifth roots of unity. Taking $m = o(\epsilon)$, we deduce that $c(\alpha) = \lcm[5,6,m]$, so in particular, each term of $\alpha$ lies in $\Q_{c(\alpha)}$ and $G = \Gal(\Q_{c(\alpha)}/\Q(\alpha))$.
	
	Set $n = c(\alpha)$ and $\sigma(\zeta_n) = \zeta_n^x$ for some $x \in \Z_n^\times$. We have deduced that $\alpha = \epsilon(\zeta_6^c + \zeta_5^a + \zeta_5^b)$ where $c \in \{1,5\}$ and $a,b \in \{1,2,3,4\}$ for $a \neq b$. We wish to obtain more information about $a,b,c,$ and $x$. We have\\
	
	$$
	\alpha - \sigma(\alpha) = \epsilon(\zeta_6^c + \zeta_5^a + \zeta_5^b) + \epsilon^\sigma(-(\zeta_6^c)^\sigma - (\zeta_5^a)^\sigma - (\zeta_5^b)^\sigma)
	$$ $$
	= \epsilon\big[\zeta_6^c + \zeta_5^a + \zeta_5^b + \epsilon^{x-1}(-(\zeta_6^c)^\sigma - (\zeta_5^a)^\sigma - (\zeta_5^b)^\sigma)\big].
	$$
	In order for $\alpha - \sigma(\alpha)$ to have the form in \eqref{E:4}, we must have $\epsilon^{x-1} = -1$, $(\zeta_6^c)^\sigma = \zeta_6^{5c}$, and $$\{\zeta_5^a , \zeta_5^b\} \cap \{(\zeta_5^a)^\sigma, (\zeta_5^b)^\sigma\} = \emptyset.$$ By considering how $\Gal(\Q_5)$ acts on the primitive fifth roots of unity, one sees that the only possible choice of $a$ and $b$ is given by $\{a,b\} = \{1,4\}$ or $\{2,3\}$. Since $2 \cdot 4 \equiv 3 \m 5$ and $3 \cdot 4 \equiv 2 \m 5$, we see that the two fifth roots of unity in $\alpha/\epsilon$ are fourth powers of each other, that is, $b = 4a$.
	
	=We have reduced to the case where $\alpha = \epsilon(\zeta_6^c + \zeta_5^a + \zeta_5^{4a})$ and $\sigma(\alpha) = -\epsilon(\zeta_6^{5c} + \zeta_5^{2a} + \zeta_5^{3a})$. We have exactly two possibilities for how $\sigma$ acts on $\Gal(\Q_5)$:
	
	\begin{itemize}
		\item[(i)] $\zeta_5^\sigma = \zeta_5^2$ and $(\zeta_5^4)^\sigma = \zeta_5^3$, or
		\item[(ii)] $\zeta_5^\sigma = \zeta_5^3$ and $(\zeta_5^4)^\sigma = \zeta_5^2$.
	\end{itemize}
	
	If we are in situation (i), then by the Chinese Remainder Theorem, we conclude that the restriction of $\sigma$ to $\Q_{30}$ is given by $\zeta_{30}^\sigma = \zeta_{30}^{17}$. We wish to describe the possible choices of $\epsilon$ for which such a $\sigma$ exists. In particular,  $x$ must satisfy
	
	\begin{equation}\label{E:5}
		x \equiv 17 \m 30,
	\end{equation}
	and since $\epsilon^{x-1} = -1$, we obtain 
	\begin{equation}\label{E:6}
		x \equiv \frac{m}{2}+1 \m m.
	\end{equation}
	We can combine \eqref{E:5} and \eqref{E:6} to obtain a necessary condition on $m$. In particular, $m$ must be of the form
	
	\begin{equation}\label{E:7}
		m = \dfrac{32 + 60k}{2\ell - 1},
	\end{equation}
	where $k$ and $\ell$ range over all integers which make $m$ positive.

	It is clear from \eqref{E:7} that $m$ is divisible by 4 but not divisible by 3 and 5. To see that the converse holds, we set $m' = m/4$ and consider a rearrangement of \eqref{E:7} after factoring out 4:
	
	\begin{equation}\label{E:8}
		8 + m' = 2m'\ell - 15k.
	\end{equation}
	If $m'$ is divisible by either 3 or 5, then $\gcd(2m' , 15) \neq 1$, and so, viewing \eqref{E:8} as a linear Diophantine equation, an integer solution $(\ell, k)$ cannot exist. 
	
	Situation (ii) is entirely similar to the first, the difference being that $\sigma$ acts on $\Q_{30}$ as $\zeta_{30}^\sigma = \zeta_{30}^{23}$ by another application of the Chinese Remainder Theorem. One sees that this automorphism is the inverse of the $\sigma$ considered in situation (i), so there is no need to investigate this situation separately.
	
	The claim that $|G| = 4$ follows from the fact that $\sigma$ and $\sigma^{-1}$ are the unique automorphisms of $G$ not arising from a permutation of the terms of $\alpha$. On the other hand, there can only be one automorphism that permutes the terms, namely the one that swaps $\epsilon\zeta_5$ and $\epsilon\zeta_5^4$. However, this automorphism is given by $\sigma^2$. $\square$ 
	
	\section{Cyclotomic Integers Of Length 4}
	
	A counterexample to Conjecture \ref{C:3} of length 4 is given by
	$$\alpha = i(\zeta_6 + \zeta_7 + \zeta_7^2 + \zeta_7^4).$$
	Here, $c(\alpha) = 84$, and there exists $\sigma \in \Gal(\Q_{84}/\Q(\alpha))$ of order 6 given by $\sigma(\zeta_{84}) = \zeta_{84}^{59}$. All counterexamples of length 4 will take a similar form. However, the case analysis is complicated by the fact that there is not one but rather three minimal vanishing sums having weight 8. These minimal vanishing sums have been classified by Conway and Jones in \cite{ConJon}. In particular, a minimal vanishing sum of weight 8 will take one of the following forms:
	
	\begin{equation}\label{E:9}
		\epsilon(\zeta_6 + \zeta_6^{5} + \sum_{i=1}^6 \zeta_7^{i}) = 0,
	\end{equation}
	\begin{equation}\label{E:10}
		\epsilon\big[(\zeta_6 + \zeta_6^5) + \zeta_5(\zeta_6 + \zeta_6^5) + \zeta_5^2(\zeta_6 + \zeta_6^5) + \zeta_5^3 + \zeta_5^4)\big] = 0,
	\end{equation}
	or
	\begin{equation}\label{E:11}
		\epsilon\big[(\zeta_6 + \zeta_6^5) + \zeta_5^2(\zeta_6 + \zeta_6^5) + \zeta_5^3(\zeta_6 + \zeta_6^5) + \zeta_5 + \zeta_5^4)\big] = 0,
	\end{equation}
	where $\epsilon$ is some root of unity.
	
	Let $\alpha = \epsilon_1 + \epsilon_2 + \epsilon_3 + \epsilon_4$ be a cyclotomic integer of length 4, and set $G = \Gal(\Q(\epsilon_1, \epsilon_2,\epsilon_3, \epsilon_4)/\Q(\alpha))$. For $\sigma \in G$, we may partition $\alpha - \sigma(\alpha)$ into minimal vanishing sums. If for each $\sigma \in G$, we may partition the sum into four sums of weight 2, then $G$ permutes the terms of $\alpha$. In this case, $G$ is an abelian subgroup of $S_4$ and has order at most 4. If we can produce a partition of the terms including a set of size 3, then it is possible to replace $\alpha$ with a sum of weight 3 or less as was the case in Theorem \ref{T:4}, a contradiction to the assumption $l(\alpha) = 4$. Such a partition cannot contain a set of size 4 either since no minimal vanishing sums of weight 4 exist.

	A slightly more complicated possibility is that $\alpha - \sigma(\alpha)$ may be partitioned into a set of size 6 and a set of size 2. In this case, Theorem \ref{T:4} implies that $\alpha$ must take the form
	
	$$
	\alpha = \epsilon(\zeta_6^x + \zeta_5^y + \zeta_5^{4y}) + \delta
	$$
	for $x \in \Z_6^\times$, $y \in \Z_5^\times$, and roots of unity $\epsilon$ and $\delta$. One verifies that
	
	$$
	\sigma(\alpha) = -\epsilon(\zeta_6^{5x} + \zeta_5^{2y} + \zeta_5^{3y}) + \delta.
	$$
	We claim that $\delta^\sigma = \delta$, for otherwise, we would have that $\delta$ is a product of $\epsilon$ times either a primitive fifth or sixth root of unity. In either case, this contradicts the assumption that $l(\alpha) = 4$. Theorem 4 now implies that $|G| = 4$.
	
	\begin{theorem}\label{T:5}
		Let $\alpha$ be a cyclotomic integer of length 4. Then, $$[\Q_{c(\alpha)}:\Q(\alpha)] \leq 6.$$ Moreover, equality holds iff $\alpha$ is of the form
		
		$$
		\epsilon(\zeta_6^c + \zeta_7^a  + \zeta_7^{2a} +  \zeta_7^{4a}),
		$$
		where $c \in \Z_6^\times$, $a \in \Z_7^\times$, and $\epsilon$ is a root of unity whose order is a multiple of 4 not divisible by either 3 or 7. 
		
	\end{theorem}
	
	\noindent\emph{Proof:} We have shown above that the only possible counterexamples to the conjecture having weight 4 are exactly the $\alpha$ for which there exists $\sigma \in G$ such that $\alpha - \sigma(\alpha)$ is itself a minimal vanishing sum of weight 8. We claim that the only counterexamples arise from \eqref{E:9}. Handling case \eqref{E:9} is entirely analogous to the work shown in the proof of Theorem \ref{T:4}, so we are only concerned with ruling out the existence of counterexamples arising from \eqref{E:10} or \eqref{E:11}. Moreover, our proof that no counterexamples arise from \eqref{E:10} will be practically identical to the proof that none arise from \eqref{E:11}, so we concern ourselves with \eqref{E:10} only.
	
	Let $\alpha = \epsilon_1 + \epsilon_2 + \epsilon_3 + \epsilon_4$ be a cyclotomic integer of length 4, and set $G = \Gal(\Q(\epsilon_1, \epsilon_2,\epsilon_3, \epsilon_4)/\Q(\alpha))$. For $\sigma \in G$, assume that $\alpha - \sigma(\alpha)$ is of the form in \eqref{E:10}. We consider the various ways to choose $\alpha$. Our assumption that $l(\alpha) = 4$ forces $\alpha/\epsilon$ to have exactly one term which is a primitive sixth root of unity. We assume this term is $\zeta_6^c$ for some $c \in \Z_6^\times$. Also, $\alpha/\epsilon$ must have exactly one term which is a primitive fifth root of unity. Otherwise, $\sigma(\alpha)$ would have three terms that are multiples of either $\zeta_6^c$ or $\zeta_6^{5c}$ which would allow us to replace $\alpha$ with a sum of three roots of unity. For instance, if $\alpha = \epsilon(\zeta_6 + \zeta_5^{3} + \zeta_5^{4} + \zeta_6\zeta_5)$, then
	
	$$
	\sigma(\alpha) = -\epsilon(\zeta_6^5 + \zeta_6^5\zeta_5 + \zeta_6^5\zeta_5^2 + \zeta_6\zeta_5^2) = -\epsilon\big(\zeta_6^5(1 + \zeta_5 + \zeta_5^2) + \zeta_6\zeta_5^2\big)
	$$ $$
	= -\epsilon\big(\zeta_6^5(-\zeta_5^3 - \zeta_5^4) + \zeta_6\zeta_5^2\big),
	$$
	contradicting our assumption that $l(\alpha) = 4$. We may assume that this term is $\zeta_5^b$ for some $b \in \Z_5^\times$. Thus, $\alpha/\epsilon$ must be a sum consisting of a primitive sixth root of unity, a primitive fifth root of unity, and two primitive thirtieth roots of unity. Moreover, these two primitive thirtieth roots of unity will be of the form $\zeta_6^{a_1}\zeta_5$ and $\zeta_6^{a_2}\zeta_5^2$ for $a_1, a_2 \in \Z_6^\times$. Otherwise, $\alpha$ would again have length less than 4. We may now deduce
	
	$$
	\alpha = \epsilon(\zeta_6^c + \zeta_6^{a_1}\zeta_5 + \zeta_6^{a_2}\zeta_5^2 + \zeta_5^b)
	$$
	and
	$$
	\sigma(\alpha) = -\epsilon(\zeta_6^{5c} + \zeta_6^{5a_1}\zeta_5^{y} + \zeta_6^{5a_2}\zeta_5^{2y} + \zeta_5^{by})
	$$
	for $c, a_1, a_2 \in \Z_6^\times$, $y \in \Z_5^\times$, and $b, by \in \{3,4\}$.
	
	We claim that $\sigma$ has order at most 2 and that $G = \la \sigma \ra$. First, no element of $G$ permutes the terms of $\alpha$ nontrivially. Indeed, the only possible permutation would swap the two primitive thirtieth roots of unity. However, this would require a $\tau \in G$ to swap $\zeta_5$ and $\zeta_5^2$, an impossibility. By comparing the terms of $\alpha$ to the minimal vanishing sum from Theorem \ref{T:4}, it is also impossible for $\alpha - \tau(\alpha)$ to have a subsum which is a minimal vanishing sum of weight 6. We now have that $\alpha - \tau(\alpha)$ is a minimal vanishing sum for all nonidentity $\tau$. However, $\alpha - \tau(\alpha)$ is never of the form in \eqref{E:9} or \eqref{E:11}, for otherwise, there would then exist a root of unity $\epsilon'$ such that $\epsilon'\alpha$ is a subsum of either \eqref{E:9} or \eqref{E:11}, but it is straightforward to verify that no such $\epsilon'$ can exist. It now follows that $\sigma^2 = 1$, and no other nontrivial automorphisms exist. So, $|G| \leq 2$. $\square$.\\
	
	\section{The Asymptotic Bound}
	
	%\begin{lemma}
	%	The largest possible order of an abelian subgroup of $S_n$ is either $2\cdot3^{\lfloor \frac{n}{3} \rfloor}$ or $3^{\lfloor \frac{n}{3} \rfloor}$ depending on whether $n \equiv 2 \m 3$ or not.
	%\end{lemma}
	
	%\paragraph{} \textit{Proof:} Let $A$ be an abelian subgroup of $S_n$ having orbits $\mathcal{O}_i$ for $i=1,\ldots,\ell$. Since abelian permutation groups act semiregularly, we have $|A| \leq \prod_{i=1}^\ell |\mathcal{O}_i|$. It then suffices to assume $A$ takes the form $\prod_{i=1}^\ell \la g_i \ra$ where the $g_i$ are mutually disjoint cycles of length $|\mathcal{O}_i|$. If any $g_i$ in the product is a cycle of length at least 5, then we may obtain a larger subgroup by replacing this cycle with two cycles of length $\lfloor \frac{m}{2} \rfloor$ and of length $ m - \lfloor \frac{m}{2} \rfloor$ since the product of these two numbers will then exceed $m$. It now follows that if $A$ has maximal order in $S_n$, then each $g_i$ is either a 2, 3, or 4-cycle. If any $g_i$ is a 4-cycle, then replacing the $g_i$ with two 2-cycles will not affect the order of $A$. Since $2^{n/2} < 3^{n/3}$ for all $n > 0$, we can replace all but at most one 2-cycles with an appropriate number of 3-cycles. Computing the order of $A$ gives the lemma. $\square$ \\
	
	We will now give the proof of the main theorem. We define any two sums of roots of unity $\alpha_1$ and $\alpha_2$ to be equivalent up to a rotation if there exists a root of unity $\epsilon$ such that $\epsilon\alpha_1 = \alpha_2$. Otherwise, we will call $\alpha_1$ and $\alpha_2$ distinct. The below lemma bounds the number of minimal vanishing sums of a given weight up to a rotation.
	
	\begin{lemma}\label{L:3}
		The number of minimal vanishing sums up to a rotation having weight $k$ is at most $e^{ck^2}$ for some $c \in \R$.
	\end{lemma}
	
	\noindent\emph{Proof:} Let $\alpha$ be an arbitrary minimal vanishing sum of weight $k$. By Theorem 1 in \cite{Mann}, there exists a root of unity $\epsilon$ such that the terms of $\epsilon\alpha$ all have square-free order (see also \cite{LamLe} for a ring-theoretic proof). Further, we can rotate $\alpha$ by a root of unity so that the largest prime dividing the order of any term is at most $k$. To see this, assume that 
	
	$$
	\alpha = \sum_{i=1}^k \epsilon_i,
	$$
	where each $\epsilon_i$ has square-free order. Now, let $p$ be the largest prime dividing the order of any term. By factoring the $p$-part out of each term, we may rewrite $\alpha$ as follows:
	
	\begin{equation}\label{E:12}
		\alpha = \sum_{i = 0}^{p-1} f_i \zeta_p^i,
	\end{equation}
	where each $f_i$ is a sum of roots of unity having order not divisible by $p$. By the linear independence of the primitive $p$th roots of unity, it follows that $f_0 = f_1 = \ldots = f_{p-1}$. If the $f_i$'s are zero, minimality implies that exactly one of the $f_i$'s is not identically zero. If this occurs at $f_j$, then $\alpha\zeta_p^{-j}$ is a rotationally equivalent minimal vanishing sum with terms having order divisible only by strictly smaller primes. Inductively, we may now assume that each of the $f_i$'s in \eqref{E:12} is nonzero so that $w(f_i) \geq 1$ for each $i$. Now,
	
	$$
	k = w\bigg(\sum_{i = 0}^{p-1} f_i \zeta_p^i\bigg) = \sum_{i=0}^{p-1} w(f_i) \geq p,
	$$
	which verifies the claim.
	
	Set $\ell = \prod_{i=1}^{\pi(k)} p_i$. We have shown in the previous paragraph that any minimal vanishing sum of weight $k$ is rotationally equivalent to a sum of $\ell$th roots of unity. The number of ways of choosing a sum of $k$ $\ell$th roots of unity is seen to be $\binom{\ell + k - 1}{k}$. Now, $\ell = e^{\vartheta(k)}$ where $\vartheta(k) = \sum_{p \leq k} \log(p)$ is the Chebyshev function. It is shown in \cite{RossSchoen} that $\vartheta(k) < 1.000028k$ for all $k \geq 1$. Set $c = 1.000028$. We have
	
	$$
	\binom{\ell + k - 1}{k} = \dfrac{\prod_{i=0}^{k-1} (\ell + i)}{k!} < \dfrac{\prod_{i=0}^{k-1} (e^{ck} + i)}{k!} \leq e^{ck^2},
	$$
	which proves the lemma. $\square$ \\
	
	\begin{xrem} 
		Mann proved in Theorem 2 of \cite{Mann} that the number of minimal vanishing sums of a given weight is finite using this same approach. However, he did not attempt to produce an upper bound for the number of such minimal vanishing sums.
	\end{xrem}
	
	We are now able to produce our desired bound.\\
	
	\noindent\emph{Proof of Main Theorem:} Let $\alpha$ be a cyclotomic integer of length $k$, and set $\alpha = \sum_{i=1}^k \epsilon_i$ for roots of unity $\epsilon_i$. Define $K$ to be the extension of $\Q$ generated by all of the $\epsilon_i$, and set $G = \Gal(K/\Q(\alpha))$. Set $X = \{\epsilon_i\}_{i=1}^k$. We view $X$ as a multiset since some terms of $\alpha$ may be repeated. For each $\sigma \in G$, we define the multiset $\sigma(X) = \{\sigma(\epsilon_i) \; | \; \epsilon_i \in X\}$. Our approach will begin with bounding the possible choices of $\sigma(X)$. We proceed in steps.
	
	\emph{Step 1:} Choose a $\sigma \in G$. We may partition the disjoint union of multisets $X \sqcup \sigma(X)$ into terms which form minimal vanishing subsums of $\alpha - \sigma(\alpha)$, say $X \sqcup \sigma(X) = Y_1 \sqcup \ldots \sqcup Y_\ell$. Similarly to the proof of Lemma \ref{L:2}, our length assumption forces each $Y_i$ to consist of the same number of elements of $X$ as of $\sigma(X)$, i.e., $|X \cap Y_i| = |\sigma(X) \cap Y_i|$ for each $i$. %We wish to bound the possible choices of $Y_i - (X \cap Y_i)$.
	
	\emph{Step 2:}  Choose a partition $X = X_1 \sqcup \ldots \sqcup X_\ell$. For a given $X_i$, we wish to count the possible choices of $\tilde{X}_i$ such that $|X_i| = |\tilde{X_i}|$ and the sum over $X_i \sqcup \tilde{X}_i$ is a minimal vanishing sum. The $\tilde{X}_i$'s will play the role of a possible $\sigma(X) \cap Y_i$. Let $\alpha_i$ denote the sum over the elements of $X_i$, and set $k_i = |X_i|$. By Lemma \ref{L:3}, there are less than $e^{4ck_i^2}$ distinct minimal vanishing sums of weight $2k_i$. For a fixed minimal vanishing sum $\beta_i$ of weight $2k_i$, we claim that $\alpha_i$ can be a subsum of at most $2k_i$ rotations of $\beta_i$. To see this, assume that $\alpha_i$ is a subsum of $\epsilon\beta_i$ for some root of unity $\epsilon$, so that $\epsilon^{-1}\alpha_i$ is a subsum of $\beta_i$. 
	Any rotation of $\alpha_i$ is determined by where the rotation sends one of the terms, and there are at most $2k_i$ rotations of a given root of unity to some term of $\beta_i$. There then can be at most this many rotations of $\alpha_i$ to some subsum of $\beta_i$. It now follows that there are at most $2k_i$ choices of $\epsilon$, which proves the claim.
	Thus, $\alpha_i$ can be a subsum of at most $2k_i e^{4ck_i^2}$ minimal vanishing sums of weight $2k_i$. It follows that there are at most this many choices for the multiset $\tilde{X}_i$.
	
	\emph{Step 3:} Accounting for each $i$ simultaneously, we see that the number of possibilities for $\sigma(X)$ which respect the chosen partition is at most
	
	\begin{equation}\label{E:13}
		\prod_{i=1}^\ell 2k_i e^{4ck_i^2} \leq 3^k \cdot e^{4ck^2}.
	\end{equation}
	Here, we used the fact that $\sum_i 2k_i = 2k$ implies that $\prod_i 2k_i \leq 3^{\lceil \frac{2k}{3}\rceil} \leq 3^k$.
	
	\emph{Step 4:} The bound in \eqref{E:13} was dependent on the choice of a partition for $X$. As $|X| = k$, there are $B_k$ choices for a partition of $X$, where $B_k$ is the $k$th Bell number. Thus, using \eqref{E:13}, we see that the total number of possibilities for $\sigma(X)$ is at most
	
	\begin{equation}\label{E:14}
		B_k\bigg(3e^{4ck}\bigg)^k.
	\end{equation} It has been shown in \cite{BerTass} that $B_k < \big(\frac{0.792k}{\log(k+1)}\big)^k$
	for all choices of $k$, so we may substitute this expression in for $B_k$ above. 
	
	\emph{Step 5:} Now that we have determined a bound for the number of possibilities for $\sigma(X)$, we can now state the bound for $|G|$ in terms of $k$. For a given choice of $\sigma(X)$, the number of different ways for an automorphism to assign terms of $X$ to $\sigma(X)$ is at most the number of bijections between these multisets, which may be as large as $k!$. Thus, multiplying the expression in \eqref{E:14} by $k!$ gives an upper bound to the number of ways an automorphism of $G$ can move the terms of $\alpha$. Since each automorphism of $G$ is determined by how the terms of $\alpha$ are moved, this is a bound for $|G|$. Since $|G| \geq |\Gal(\Q_{c(\alpha)}/\Q(\alpha))|$, the theorem now follows. $\square$ \\

	\section{Conclusion}
	
	We have seen in Theorems \ref{T:3}, \ref{T:4}, and \ref{T:5} that counterexamples to the conjecture having length at most 4 are highly restricted. This leads us to wonder whether these theorems can be used to produce a second, possibly more elementary proof to the special cases of Conjecture \ref{C:2} where $\chi(1) \leq 3$. Such a proof could likely be modified to handle cases involving slightly larger character degrees as well. However, we are unable to produce any simpler proof utilizing our results at the moment.
	
	We remark that the bound on $d(k)$ in main theorem seems very crude, and we suspect that this bound can be improved drastically. One potential approach may be to tighten the bound obtained in Lemma \ref{L:3}. To this date, minimal vanishing sums of weight at most 21 have been classified in \cite{ChrisDyk}, \cite{ConJon}, \cite{Mann}, and \cite{PoonRub}. However, we are unaware of any results about the asymptotic behavior on the number of distinct minimal vanishing sums. The results from \cite{PoonRub} seem to indicate that such a bound might be $O(e^{ak})$ for some $a > 0$ as opposed to $O(e^{ck^2})$ as we have shown in Lemma \ref{L:3}.  Any results in this direction would be interesting in their own right as minimal vanishing sums seem to have a wide range of applications.

\end{document}